\documentclass[12pt]{amsart}
\usepackage{amscd,amsmath,amssymb,amsfonts}
\usepackage[all]{xy}
\theoremstyle{plain}
\newtheorem{thm}{Theorem}
\newtheorem{lem}[thm]{Lemma}

\newtheorem{prop}[thm]{Proposition}

\theoremstyle{definition}
\newtheorem{defn}[thm]{Definition}

\newtheorem{rmk}[thm]{Remark}
\newtheorem{rmks}[thm]{Remarks}

\numberwithin{thm}{section}
\numberwithin{equation}{section}

\newcommand{\ml}[2]{\begin{multline}\label{#1}#2 \end{multline}}
\newcommand{\ga}[2]{\begin{gather}\label{#1}#2 \end{gather}}

\newcommand{\sA}{{\mathcal A}}
\newcommand{\sB}{{\mathcal B}}
\newcommand{\sC}{{\mathcal C}}

\newcommand{\sE}{{\mathcal E}}

\newcommand{\sK}{{\mathcal K}}

\newcommand{\sO}{{\mathcal O}}

\renewcommand{\H}{{\mathbb H}}

\begin{document}

\title[Partial connection for $p$-torsion line bundles]{Partial connection for $p$-torsion line bundles in characteristic $p>0$}
\author{H\'el\`ene Esnault}
\address{Universit\"at Duisburg-Essen, Mathematik, 45117 Essen, Germany}
\email{esnault@uni-due.de}
\date{August 7, 2006}
\dedicatory{To S. S. Chern, in memoriam}
\thanks{Partially supported by  the DFG Leibniz Preis}
\begin{abstract}
The aim of this brief note is to give a construction for $p$-torsion line bundles in characteristic $p>0$ which plays a similar r\^ole as the standard connection on an $n$-torsion line bundle in characteristic 0. 
\end{abstract}
\maketitle
\begin{quote}

\end{quote}

\section{Introduction}
In \cite{E1} (see also \cite{E2}) we gave an algebraic construction of characteristic classes of vector bundles with a flat  connection $(E,\nabla)$ on a smooth algebraic variety $X$ defined over a field $k$ of characteristic 0.  Their value at the generic point ${\rm Spec}(k(X))$   was studied and redefined in 
\cite{BE1}, and then applied in \cite{BE2} to establish a Riemann-Roch formula. 
One way to understand Chern classes of vector bundles (without connection) is via the Grothendieck splitting principle: if the receiving groups $\oplus_n H^{2n}(X,n)$ of the classes form a cohomology theory which is a ring and is functorial in $X$, then via the Whitney product formula it is enough to define the first Chern class. Indeed, on the flag bundle $\pi: {\rm Flag}(E)\to X$, $\pi^*(E)$ acquires a complete flag $E_i\subset E_{i+1}\subset \pi^*(E)$ with $E_{i+1}/E_i$ a line bundle, and $\pi^*: H^{2n}(X,n) \to H^{2n}({\rm Flag}(E),n)$ is injective, so it is enough to construct the classes on ${\rm Flag}(E)$.  However, if $\nabla$ is a connection on $E$, $\pi^*(\nabla)$ does not stabilize the flag $E_i$. So the point of \cite{E1} is to show that there is a  differential graded algebra $A^\bullet$ on ${\rm Flag}(E)$, together with 
a morphism of differential graded algebras $\Omega^\bullet_{{\rm Flag}(E)}\xrightarrow{\tau} A^\bullet$, so that $R\pi_*A^\bullet \cong \Omega^\bullet_X$ and so that the operator defined by the composition
$\pi^*(E) \xrightarrow{\pi^*(\nabla)}  \Omega^1_{{\rm Flag}(E)} \otimes_{\sO_{{\rm Flag}(E)}} \pi^*(E)
\xrightarrow{\tau \otimes 1} A^1\otimes_{\sO_{{\rm Flag}(E)}} \pi^*(E)$ stabilizes $E_i$.  We call the induced operator $\nabla_i: E_i \to A^1\otimes_{\sO_{{\rm Flag}(E)}} E_i$ a ({\it flat})  $\tau$-{\it connection}. So it is a $k$-linear map which fulfills the $\tau$-{\it Leibniz formula} 
\ga{1.1}{\nabla_i(\lambda \otimes e)= \tau d(\lambda)\otimes e + \lambda \nabla_i(e) } for $\lambda$ a local section of $\sO_{{\rm Flag}(E)}$ and $e$ a local section of $E_i$. It is flat when $0=\nabla_i\circ \nabla_i\in H^0(X, A^2\otimes_{\sO_X}\sE nd(E))$, with the appropriate standard sign for the derivation of forms with values in $E_i$. The last point is then to find the correct cohomology which does not get lost under $\pi^*$. It is a generalization of the classically defined group 
\ga{1.2}{\H^1(X, \sO^\times_X\xrightarrow{d {\rm log}} \Omega^1_X\xrightarrow{d}\Omega^2_X\xrightarrow{d} \cdots)} of isomorphism classes of rank one line bundles on $X$ with a flat connection. 

A typical example of such a connection is provided by a torsion line bundle: if $L$ is a line bundle on $X$ which is $n$-torsion, that is which is endowed with an isomorphism $L^n\cong \sO_X$, then the isomorphism yields an $\sO_X$-\'etale algebra structure on $\sA=\oplus_0^{n-1} L^i$, hence a finite  \'etale covering $\sigma: Y={\rm Spec}_{\sO_X} \sA \to X$, which is a principal bundle under the group scheme $\mu_n$ of $n$-th roots of unity, thus is Galois cyclic as soon as  $\mu_n\subset k^\times$. Since the $\mu_n$-action commutes with the differential $d_Y: \sO_Y\to \Omega^1_Y=\sigma^*\Omega^1_X$, it defines a flat connection $\nabla_L: L\to \Omega^1_X\otimes_{\sO_X} L$. Concretely, if $g_{\alpha, \beta}\in \sO_X^\times$ are local algebraic transition functions for $L$, with trivialization
\ga{1.3}{g^n_{\alpha, \beta}=u_\beta u_\alpha^{-1}, u_\alpha\in \sO^\times_X,}
then 
\ga{1.4}{(g_{\alpha, \beta}, \frac{1}{n} \frac{du_\alpha}{u_\alpha})\in \big(\sC^1(\sO^\times_X)\times \sC^0((\Omega^1_X)_{{\rm clsd}})\big)_{{d {\rm log}-\delta}}\\
\frac{dg}{g}=\delta(\frac{du}{u})\notag}
is a Cech cocyle for the class 
\ga{1.5}{(L,\nabla_L)\in \H^1(X, \sO^\times_X\xrightarrow{d {\rm log}} \Omega^1_X\xrightarrow{d}\Omega^2_X\xrightarrow{d} \cdots).} 
Clearly \eqref{1.4} is meaningless if the characteristic $p$ of $k$ is positive and divides $n$. The purpose of this short note is to give an Ersatz of this canonical construction in the spirit of the $\tau$-connections explained above when $p$ divides $n$. 
\section{A partial connection for $p$-torsion line bundles}
Let $X$ be a scheme of finite type over a perfect field $k$ of characteristic $p>0$. 
Let $L$ be a $n$-torsion line bundle on $X$, thus endowed with an isomorphism
\ga{2.1}{\theta: L^n\cong \sO_X.}
Then $\theta$ defines an $\sO_X$-algebra structure on 
$\sA=\oplus_0^{n-1} L^i$
which is \'etale if and only if $(p,n)=1$. It defines the principal $\mu_n$-covering
\ga{2.2}{\sigma: Y={\rm Spec}_{\sO_X} \sA\to X}
which is \'etale if and only if $(p,n)=1$, else decomposes into 
\ga{2.3}{\sigma: Y\xrightarrow{\iota} Z\xrightarrow{\sigma'} X}
with $\sigma'$ \'etale and  $\iota$ purely inseparable. More precisely, if $n=m\cdot p^r, (m,p)=1$, and $M=L^{p^r}$, $\theta$ defines an $\sO_X$-\'etale algebra structure on
$\sB=\oplus_0^{m-1} M^i$, which defines $\sigma': Z={\rm Spec}_{\sO_X} \sB\to X$ as an (\'etale) $\mu_m$-principal bundle. The isomorphism $\theta$ also defines an isomorphism  $(L')^{p^r}\cong \sO_Z$ as it defines the isomorphism $(\sigma')^*(M)\cong \sO_Z$, where  $L'=(\sigma')^*(L)$. So  $\sC=\oplus_0^{p^r-1} (L')^i$
becomes a finite purely inseparable $\sO_Z$-algebra defining the principal $\mu_{p^r}$-bundle $\iota: Y={\rm Spec}_{\sO_Z} \sC \to Z$. \\ \ \\
If $(n, p)=1$, that is if $r=0$,  the formulae 
\eqref{1.3}, \eqref{1.4} define $(L,\nabla)$ as in \eqref{1.5}.
We assume from now on that $(n,p)=p$. Then,  as is well known, as a consequence of \eqref{1.3} one sees that the form 
\ga{2.4}{\omega_L:=\frac{du_\alpha}{u_\alpha} \in \Gamma(X, \Omega^1_X)_{{\rm clsd}}^{{\rm Cartier}=1}}
 is globally defined and Cartier invariant. 
Let $e_\alpha$ be local generators of $L$, with transition functions $g_{\alpha, \beta}$ with $e_\alpha=g_{\alpha, \beta} e_\beta$.
The isomorphism $\theta$ yields a trivialization
\ga{2.5}{\sigma^*L\cong \sO_Y}
thus local units $v_\alpha$ on $Y$ with 
\ga{2.6}{v_\alpha\in \sO^\times_Y, \ g_{\alpha, \beta}=v_\beta v_\alpha^{-1}\\
{\rm so \ that} \ 1=v_\alpha \sigma^*(e_\alpha)=v_\beta \sigma^*(e_\beta). \notag}
\begin{defn} \label{defn2.1}
One defines the $\sO_X$-coherent sheaf $\Omega^1_L$ as the subsheaf of $\sigma_*\Omega^1_Y$ spanned by ${\rm Im}(\Omega^1_X)$ and $\frac{dv_\alpha}{v_\alpha}$. 
\end{defn}
\begin{lem}\label{lem2.2} $\Omega^1_L$ is well defined and one has the 
 exact sequence
\ga{2.7}{0\to \sO_X\xrightarrow{\cdot \omega_L} \Omega^1_X\xrightarrow{\sigma^*} \Omega^1_L \xrightarrow{s} \sO_X\to 0\\
s(\frac{dv_\alpha}{v_\alpha})=1.\notag}

 \end{lem}
\begin{proof} The relation \eqref{2.6} implies
\ga{2.8}{\frac{dg_{\alpha, \beta}}{g_{\alpha, \beta}}=\frac{dv_\beta}{v_\beta}-\frac{dv_\alpha}{v_\alpha}\\
{\rm so} \ \frac{dv_\beta}{v_\beta} \equiv \frac{dv_\alpha}{v_\alpha} \ \in 
\sigma_*\Omega^1_Y/ {\rm Im}(\Omega^1_X).\notag}
Hence the sheaf $\Omega^1_L$ is well defined. If $e'_\alpha$ is another basis, then one has $e_\alpha=w_\alpha e'_\alpha$ for local units $w_\alpha \in \sO_X^\times$. The new $v_\alpha$ are then multiplied by local units in $\sO^\times_X$, so the surjection $s$ is well defined. 
 It remains to see that ${\rm Ker}(\sigma^*)={\rm Im}(\cdot \omega_L)$.  By definition, on the open of $X$ on which $L$ has basis $e_\alpha$, one has 
\ga{2.9}{Y={\rm Spec} \ \sO_X[v_\alpha]/(v^n_\alpha -u_\alpha).} This
implies $\Omega^1_Y=\langle {\rm Im}(\Omega^1_X), dv_\alpha\rangle_{\sO_Y}/\langle du_\alpha\rangle_{\sO_Y}$ on this open and finishes the proof.

\end{proof}
\begin{rmks} \label{rmk2.3} \begin{itemize}
\item[1)]Assume for example that $X$ is a smooth projective curve of genus $g$, and $n=p$.
Recall that $0\neq  \omega_L\in \Gamma(X, \Omega^1_X)$. In particular, if $g\ge 2$, necessarily $0\neq \Omega^1_X/\sO_X\cdot \omega_L$ is supported in codimension 1. So $\Omega^1_L$ contains a non-trivial torsion subsheaf. 
\item[2)] The sheaf $\Omega^1_L$ lies in $\sigma_*\Omega^1_Y$ but is not equal to it. Indeed,
on the smooth locus of $X$ (assuming $X$ is reduced) the torsion free quotient of $\Omega^1_L$ has rank equal to the dimension of $X$, while $\sigma_*\Omega^1_Y$ has rank $n\cdot$ dimension $(X)$ on the \'etale locus of $\sigma$ (which is non-empty if $L$ itself is not a $p$-power line bundle).
\item[3)] The class in ${\rm Ext}^2_{\sO_X}(\sO_X, \sO_X)=H^2(X, \sO_X)$ defined by \eqref{2.7} vanishs. Indeed, let us decompose \eqref{2.7} as an extension of $\sO_X$ by $\Omega^1_X/\sO_X\cdot \omega_L$, followed by an extension of $\Omega^1_X/\sO_X\cdot \omega_L$ by $\sO_X\cdot \omega_L$. The first extension class in $H^1(X, \Omega^1_X/\sO_X\cdot \omega_L)$ has cocycle $\frac{dv_\beta}{v_\beta} - \frac{dv_\alpha}{v_\alpha}=\frac{dg_{\alpha, \beta}}{g_{\alpha, \beta}}$ (see \eqref{2.8}), thus is the image of the Atiyah class of $L$ in $H^1(X, \Omega^1_X)$. Thus  the second boundary to $H^2(X, \sO_X)$ dies. 
\end{itemize}
\end{rmks}
\begin{defn} \label{defn2.4}
We set $\Omega^0_L:=\sO_X$ and for $i\ge 1$ we define the $\sO_X$-coherent sheaf
$\Omega^i_L$  as the subsheaf of $\sigma_*\Omega^i_Y$ spanned by ${\rm Im}(\Omega^i_X)$ and 
$\frac{dv_\alpha}{v_\alpha}\wedge {\rm Im}(\Omega^{i-1}_X)$. 
\end{defn}
\begin{prop} \label{prop2.5}
The sheaf $\Omega^i_L$ is well defined. One has an exact sequence
\ga{2.10}{0\to \omega_L\wedge \Omega^{i-1}_X\to \Omega^i_X\xrightarrow{\sigma^*} \Omega^i_L \xrightarrow{s} \Omega_X^{i-1}\to 0\\
s(\frac{dv_\alpha}{v_\alpha}\wedge \beta)=\beta.\notag}
Furthermore, the differential $\sigma_*(d_Y)$ on $\sigma_*\Omega^\bullet_Y$ induces on $\oplus_{i\ge 0} \Omega^i_L$ the structure of a differential graded algebra $(\Omega^\bullet_L, d_L)$ so that $\sigma^*: (\Omega^\bullet_X, d_X)\to (\Omega^\bullet_L, d_L)$ is a morphism of differential graded algebras. 
\end{prop}
\begin{proof}
One proves \eqref{2.10} as one does \eqref{2.7}. One has to see that $\sigma_*(d_Y)$ stabilizes $\Omega^\bullet_L$.  As $0=d_X(\omega_L)\in \Omega^2_X, \ 0=d_Y(\frac{dv_\alpha}{v_\alpha})\in \sigma_*\Omega^2_Y$,  \eqref{2.10} extends to an exact sequence of complexes
\ml{2.11}{0\to (\omega_L\wedge \Omega^{\bullet-1}_X, -1\wedge d_X)\to (\Omega^\bullet_X, d_X)\xrightarrow{\sigma^*} \\(\Omega^\bullet_L, d_L) \xrightarrow{s} (\Omega_X^{\bullet-1}, -d_X)\to 0.
}
This finishes the proof. 

\end{proof}
\begin{rmk} \label{rmk2.6} As $\frac{dg_{\alpha, \beta}}{g_{\alpha, \beta}}\in (\Omega^1_X)_{{\rm clsd}}$ the same proof as in Remark \ref{rmk2.3}, 3) shows that the extension class ${\rm Ext}^2(\Omega^{\bullet -1}_X, \omega_L\wedge \Omega^{\bullet-1}_X)$ defined by \eqref{2.11} dies. 
\end{rmk}
In order to tie up with the notations of the Introduction, we set
\ga{2.12}{\tau=\sigma^*: \Omega^\bullet_X\to \Omega^\bullet_L.}
\begin{prop} \label{prop2.7} The formula $\nabla(e_\alpha)=-\frac{dv_\alpha}{v_\alpha} \otimes e_\alpha \in \Omega^1_L\otimes_{\sO_X} L$ defines a flat $\tau$-connection $\nabla_L$ on $L$.  
So $(L, \nabla_L)$ is a class in $\H^1(X, \sO^\times_X\xrightarrow{\tau d{\rm log}} \Omega^1_L \xrightarrow{d_L} \Omega^2_L \xrightarrow{d_L} \cdots)$, the group of isomorphism classes of line bundles with a flat $\tau$-connection.  

\end{prop}
\begin{proof}
Formula \eqref{2.6} implies that this defines a $\tau$-connection. Flatness is obvious. A Cech cocycle for $(L,\nabla_L)$ is $(g_{\alpha, \beta}, \frac{dv_\alpha}{v_\alpha} )$.
\end{proof}
\begin{rmks} \label{rmk2.8}
\begin{itemize}
\item[1)]
The same formal definitions \ref{defn2.1} and \ref{defn2.4} of $\Omega^\bullet_L$ when $(n,p)=1$ yield $(\Omega^\bullet_L, d_L)=(\Omega^\bullet_X, d_X)$, and the flat $\tau$-connection becomes the flat connection defined in \eqref{1.4} and \eqref{1.5}. So Proposition \ref{prop2.7} is a direct genealization of it. 
\item[2)] Let $X$ be proper reduced over a perfect field $k$, irreducible in the sense that $H^0(X, \sO_X)=k$, and admitting a rational point $x\in X(k)$. A generalization of torsion line bundles to higher rank bundles is the notion of Nori finite bundles, that is bundles $E$ which are trivialized over principal bundle $\sigma: Y\to X$ under a finite flat group scheme $G$ (see \cite{N2} for the original definition and also \cite{EHS} for a study of those bundles). So for the $n$-torsion line bundles considered in this section, $G\cong \mu_n$. If the characteristic of $k$ is 0, then again $\sigma$ is \'etale, the differential $d_Y: \sO_Y\to \sigma^*\Omega^1_X=\Omega^1_Y$  commutes with the action of $G$, inducing a connection $\nabla_E: E\to \Omega^1_X\otimes_{\sO_X} E$ and characteristic classes in our groups $\H^i(X, \sK_i^m\xrightarrow{d {\rm log}} \Omega^i_X\xrightarrow{d}\Omega^{i+1}_X \cdots)$ (see \cite{E1}).  If the characteristic of $k$ is $p>0$, then $\sigma$ is \'etale if and only if $G$ is smooth (which  here means \'etale), in which case one can also construct those classes. If $G$ is not \'etale, thus contains a non-trivial local subgroupscheme, then one should construct as in Proposition \ref{prop2.5} a differential graded algebra $(\Omega^\bullet_E, d_E)$ with a map  $(\Omega^\bullet_X, d_X)\xrightarrow{\tau} (\Omega^\bullet_E, d_E)$, so that $E$ is endowed naturally with a flat $\tau$-connection $\nabla_E: E\to \Omega^1_E\otimes_{\sO_X} E$. The techniques developed in \cite{E1} should then yield classes in the 
groups $\H^i(X, \sK_i^m\xrightarrow{\tau d{\rm log}} \Omega^i_E\xrightarrow{d_E} \Omega^{i+1}_E\cdots)$. 

\end{itemize}
\end{rmks}

\bibliographystyle{plain}
\renewcommand\refname{References}

\end{document}